\documentclass[a4paper,12pt,legno]{article}

\usepackage{amsmath}
\usepackage{amsfonts}
\usepackage{amsthm}
\usepackage{mathrsfs}
\usepackage{amssymb}

\theoremstyle{remark}

\usepackage{graphicx}
\usepackage[T1]{fontenc}
\usepackage[polish]{babel}

\newcommand{\Rset}{\mathbb{R}}

\begin{document}

\title{ON A CENTRAL LIMIT THEOREM FOR
SHRUNKEN WEAKLY DEPENDENT RANDOM VARIABLES}\author{Richard C. Bradley\footnote{Department of Mathematics, 
Indiana University, Bloomington, Indiana, USA.} \ and \ Zbigniew J. Jurek\footnote{Institute of Mathematics, University of Wroc\l aw, Wroc\l aw, Poland.} \,\footnote{Research funded by Narodowe Centrum Nauki (NCN) grant no Dec2011/01/B/ST1/01257.}}

\date { }

\maketitle

{\bf Abstract.}  A central limit theorem is proved for
some strictly stationary sequences of random variables
that satisfy certain mixing conditions and are
subjected to the ``shrinking operators''
$U_r(x):=[\max\{|x|-r,0\}]\cdot x/|x|,\ r \ge 0$.
For independent, identically distributed random
variables, this result was proved earlier by
Housworth and Shao.
\hfil\break

\hfil\break


\bigskip
\bigskip

\noindent \textbf{1. Introduction}

\medskip
   In the probability theory and mathematical statistics,  many classical limit laws deal with the following sequences
$$
A_n(X_1+X_2+...+X_n)+b_n = A_nX_1+...+A_nX_n +b_n,  
\eqno (1.1)
$$
where $(X_n)$ are stochastically independent random variables or random vectors, $(A_n)$ are real numbers or linear operators and $(b_n)$ are some deterministic shifts;
cf.\  Lo\`eve [1977], Feller [1971], Jurek and Mason [1993], 
and Meerschaert and Scheffler [2001].  
Thus one could ask for limits in $(1.1)$ with some 
non-linear
deformation instead of the linear operators $A_n$.

In the middle of the 1970's, there were introduced
the non-linear shrinking operations $U_r,\, r \ge 0$, on a Hilbert space H, (in short:\ \emph{s-operations}), 
as follows:
$$
U_r(0):=0 \ \mbox{and} \ \  
U_r(x):=[\max\{||x||-r,0\}]\frac{x}{||x||}\ \ \ \mbox{for} \ x\neq 0.
$$
For these operations, the limit laws of sequences
$$
U_{r_n}(X_1)+ U_{r_n}(X_2) +\ldots+U_{r_n}(X_n) +b_n,    
\eqno (1.2)
$$
were completely characterized for infinitesimal triangular arrays; cf.\ Jurek [1977] [1981]. The limits in
(1.2) were called \emph{s-selfdecomposable and s-stable} 
for only independent $X_n$'s or independent and identically distributed random variables, respectively.
Cf.\ also Jurek [1977] [1979] [1984] [1985].  

In Jurek [1981], on page 5 it was said: \emph{The s-operations have some technical justification.
For instance, when we receive a signal X then if X is $"$small$"$, i.e., X is not greater than r, we get zero,
and in the remaining cases we get only the excess, 
i.e.\ $X-r$.}

Nowadays, in mathematical finance, note that for positive real-valued variables $X$,
$U_r(X)\equiv(X-r)^+$ (positive part)
coincides with the pay-off function of European call options; cf.\ F\"ollmer and Schied [2011].
Because of these potential applications, Marc Yor on a
couple of occasions asked one of the authors (Z.J.J.) about possible results on limits in (1.2) without the assumption
of independence.
This note will provide an answer involving a Gaussian
limit in (1.2) under certain mixing assumptions.

It is interesting to note that the class $\mathcal{U}$ of all possible limit laws in (1.2)
coincides with the totality of probability distributions of the following random integrals:
$$
\int_{(0,1]}t\,dY(t) , \ \ \ \mbox{where $Y$ is an arbitrary 
L\'evy process,} \eqno (1.3)
$$
cf.\ Jurek [1984].  From (1.3), taking $Y$ to be Brownian motion, we see that the Gaussian probability
measure can appear
as a limit in $(1.2)$;  
also cf.\ Jurek [1981], Lemma 5.2.

For real-valued random variables, Housworth and Shao [2000]
proved a central limit theorem for these
``shrinking operators'' and
stochastically independent, identically distributed random variables. In this note we prove a CLT (Theorem 2.1 below)
as in $(1.2)$ for sequences $(X_n)$ that satisfy some  
strong mixing and stationarity conditions.
\bigskip
\vfill\eject

\noindent \textbf{2. Basic notions, notations and results}

\medskip
All random variables are defined on a fixed probability space $(\Omega, {\cal F}, P)$. 
For any given $\sigma$-field ${\cal {A}} \subset {\cal {F}}$,
$L_2({\cal {A}})$ denotes the family of all 
square-integrable, 
${\cal {A}}$-measurable random variables.
By $E[X]$ we denote the expected value (if it is
defined) of a given random variable $X$; and by
$\|X\|_2$ we denote the ${\cal L}^2$-norm of a given
$X \in L_2({\cal F})$.
\medskip 
  
For any two $\sigma$-fields
${\cal A}$ and ${\cal B} \subset {\cal F}$,
define the following two measures of dependence:
$$
\alpha({\cal A},{\cal B}) :=
\sup\{
|P(A \cap B) - P(A)P(B)| : A \in {\cal{A}},  B \in {\cal{B}} \}; \eqno(2.1)$$
$$
\rho({\cal A},{\cal B}) :=
\sup\{| Corr(f,g)| : f \in L_2({\cal{A}}), g \in L_2({\cal{B}})\}
\eqno(2.2)$$
where $Corr$ stands for the correlation.
\medskip

The quantity $\rho({\cal A},{\cal B})$ is called the
``maximal correlation coefficient'' of
${\cal A}$ and ${\cal B}$.
For every pair of $\sigma$-fields
${\cal A}$ and ${\cal B}$,
one has (see e.g.\ Bradley [2007, v1, Proposition 3.11(a)])
the well known inequalities
$$ 0 \leq 4\alpha({\cal A},{\cal B})
\leq \rho({\cal A},{\cal B}) \leq 1.   \eqno (2.3) $$

Now suppose that $\textbf{X} := (X_k, k \in \bf {Z})$ is a \emph{strictly stationary} sequence of (real-valued)
random variables.
(That is, for any integers $j$ and $\ell$ and any nonnegative
integer $m$, the random vectors 
$(X_j, X_{j+1}, \dots, X_{j+m})$ and
$(X_\ell, X_{\ell+1}, \dots, X_{\ell+m})$
have the same distribution.)\ \ 
For each positive integer $n$, define the following
three dependence coefficients:
$$
\alpha(n) = \alpha(\textbf{X},n) := \alpha(
\sigma(X_k, k \leq 0), \sigma(X_k, k \geq n));
\eqno (2.4) $$
$$\rho(n) = \rho(\textbf{X},n) := \rho(
\sigma(X_k, k \leq 0), \sigma(X_k, k \geq n));
\eqno (2.5) $$
$$
\rho^*(n) = \rho^*(\textbf{X},n)  := \sup \rho(
\sigma(X_k, k \in S), \sigma(X_k, k \in T)),
\eqno (2.6)
$$
where in (2.6) the supremum is taken over all pairs
of nonempty, disjoint sets $S, T \subset {\mathbb Z}$
such that
${\rm dist}(S,T) := \min_{s \in S, t \in T}|s-t|
\geq n$.
In (2.4), (2.5), and (2.6) and below, the notation
$\sigma(\, \dots)$ refers to the $\sigma$-field of events
generated by $(\, \dots)$.
In (2.6) the two sets $S$ and $T$ can be ``interlaced'',
with each one having elements between ones in the
other set.
For each $n \in {\mathbb N}$, by (2.3)--(2.6),
$$ 0 \leq 4\alpha(n) \leq \rho(n) \leq \rho^*(n) \leq 1.
\eqno (2.7) $$ 

\medskip
The (strictly stationary) sequence 
$\textbf{X}=(X_k, k \in {\mathbb Z})$ is said to be
``strongly mixing'' (or ``$\alpha$-mixing'') if
$\alpha(n) \to 0$ as $n \to \infty$,
``$\rho$-mixing'' if $\rho(n) \to 0$
as $n \to \infty$, and
``$\rho^*$-mixing'' if $\rho^*(n) \to 0$
as $n \to \infty$.
By (2.7), $\rho$-mixing implies strong mixing, and
$\rho^*$-mixing implies $\rho$-mixing.
\medskip

The strong mixing condition was introduced by
Rosenblatt [1956], the $\rho$-mixing condition was
introduced by Kolmogorov and Rozanov [1960],
and the $\rho^*$-mixing condition was apparently
first studied by Stein [1972].
(The maximal correlation coefficient
$\rho( . \, ,.)$ in (2.2) was first studied by 
Hirschfeld [1935] in a statistical context not 
particularly connected with ``stochastic processes''.)
\medskip

Recall that for each nonnegative real number $r$, the
``shrinking operator'' $U_r: {\bf R} \to {\bf R}$ is
defined as follows: 
 $U_r(0):=0$, and for $x\neq0$
\begin{displaymath}U_r(x)\ :=\ 
[\max\{|x|-r,0\}]\frac{x}{|x|}\
=\ \left\{\begin {array}{ll}
x-r & \textrm{if $x > r$} \\
0 & \textrm{if $-r \leq x \leq r \ \  \ \ \ \qquad 
$} \\
x+r  & \textrm{if $x < -r$.}
\end{array}\right.
\eqno (2.8)
\end{displaymath}

 Now we can state the main result of this paper in terms of the all above notions:
\bigskip

{\bf Theorem 2.1.} {\sl Suppose 
$\textbf{X}=(X_k, k\in\mathbb{Z})$ is a strictly stationary sequence of random variables.
For each $r \geq 0$, define the quantity
$$ G(r) := \int_0^\infty t \cdot P(|X_0| > t+r) dt.
\eqno (2.9)$$
Suppose that
$$ \forall r\ge0,  \ \ 0<G(r)<\infty  \eqno (2.10)
$$
and that
$$ \forall \varepsilon > 0,\ \
\lim_{r \to \infty} \frac{G(r + \varepsilon)}{G(r)} = 0.
\eqno (2.11)$$

Then for each $r \geq 0$, $E[|U_r(X_0)|^2] < \infty$.

Suppose also that at least one of the following
two conditions holds: \hfil\break
\noindent (i)\ \ $\rho(1) < 1$ and
$\sum_{n=1}^\infty \rho(2^n) < \infty$; or \hfil\break
\noindent (ii)\ \ $\rho^*(1) < 1$, and
$\alpha(n) \to 0$ as $n \to \infty$.

Then there exists a sequence $(r(n), n \in {\mathbb N})$
of positive numbers satisfying
$$ r(n) \to \infty\ \ {\rm as}\ \ n \to \infty, \eqno (2.12)
$$
such that as $n \to \infty$,
$$
[U_{r(n)}(X_1) +U_{r(n)}(X_2) + \ldots +U_{r(n)}(X_n)]
-\ n \cdot m_{r(n)}\
\Rightarrow\ N(0,1),
\eqno (2.13) $$
where for each $r \geq 0$, $m_r := E[U_r(X_0)]$.}
\bigskip

   A few comments on this theorem are in order.
\medskip

   Using Fubini's Theorem (see Remark 3.2(b) in
Section 3), one can express conditions (2.10) and (2.11)  
in terms of second moments as follows
$$\forall r \geq 0,\ \ 0 <  2G(r)=\mathbb{E}[|U_r(X_0)|^2] < \infty, \eqno (2.14)
$$
and
$$ \forall \varepsilon > 0,\ \
\lim_{r \to \infty} \frac{E[|U_{r + \varepsilon}(X_0)|^2]}
{E[|U_r(X_0)|^2]} = 0.  \eqno (2.15) 
$$

   The proof of Theorem 2.1 will be given in Section 4,
after some elementary but important facts concerning
the ``shrinking operators'' $U_r(\, .\, )$ are given
in Section 3.
\medskip

   In the case where the random variables $X_k$ are
independent and identically distributed, Theorem 2.1
is due to Housworth and Shao [2000, Theorem 1].
(Recall also the earlier related work of
Jurek [1982] in connection with (1.3).)\ \   
In their result, Housworth and Shao [2000, Theorem 1] 
also gave (implicitly) some
concrete extra information on the sequence of
``shrinking parameters'' $r(n),\ n \in \mathbb{N}$.
In our context, because of the dependence between
the random variables, the same information
cannot be given.
Some partial information on the $r(n)$'s can be
seen implicitly from the proof (in Section 4) of 
Theorem 2.1.
\medskip

   In their proof of their version of Theorem 2.1
for independent, identically distributed
random variables, Housworth and Shao [2000, Theorem 1]
used (2.10) and (2.11) to in essence  
establish a Lindeberg condition --- indirectly, 
in the context of a ``codified central limit theorem'' for
independent random variables in the book of Petrov [1975].
In our different context, involving dependent random
variables, we shall adapt their arguments to
verify a Lindeberg condition directly, and then use
Lindeberg CLT's in the literature under the dependence
conditions (i) and (ii) in Theorem 2.1 respectively
to complete the proof of that theorem. 
\medskip

   In Theorem 2.1, neither weak dependence condition
(i) or (ii) implies the other; see e.g.\ 
Bradley [2007, v3, Theorem 26.8(II)].
In central limit theorems for strictly stationary
$\rho$-mixing random sequences
with finite second moments, the key role of the
``logarithmic mixing rate'' condition
$\sum_{n=1}^\infty \rho(2^n) < \infty$ in condition (i)
was established by Ibragimov [1975]; and the sharpness
of that condition in that context was shown by
counterexamples (with finite second moments and barely
slower mixing rates) such as in 
Bradley [2007, v3, Chapter 34].
In the case of strictly stationary $\rho$-mixing
sequences with (in a certain standard technical sense) ``barely infinite second moments'' such that both parts
of condition (i) hold, a central limit theorem
was proved by Bradley [1988] and extended to a 
weak invariance principle by Shao [1993].
In other central limit theory under mixing conditions 
with just finite second moments, a prominent  
role of the weak dependence condition (ii) (in 
Theorem 2.1) was established in papers such as 
Peligrad [1996] [1998] and Utev and Peligrad [2003].
\medskip

The alternative mixing assumptions (i) and (ii) in 
Theorem 2.1 provide key bounds on variances of the
partial sums
$\sum_{k=1}^n U_{r(n)}(X_k)$
in the left side of (2.13).
(Those bounds are based on Lemmas 6.1 and 6.2 in the Appendix.)\ \ 
An example described in the next theorem will show what can 
``go wrong'' in Theorem 2.1 under other, seemingly nice,
mixing conditions that fail to impose such
bounds on the variances of those partial sums.
It is a ``cancellation'' example of a type that
has long been well known in central limit theory
under dependence assumptions. 
(For an old, classic, very simple ``cancellation''
example, see e.g.\ Bradley [2007, v1, Example 1.18].)    
\hfil\break

{\bf Theorem 2.2.}\ \
{\sl Suppose $0 < \lambda < 1$.
Then there exists a strictly stationary sequence
${\bf X} := (X_k, k \in {\mathbb Z})$ of random variables
that satisfies (2.10) and (2.11)
and also has the following three properties:

\noindent (a) $\alpha(\textbf{X},1) \leq \lambda$.

\noindent (b) For some $c > 0$, 
$\rho^*(\textbf{X},n) = O(e^{-cn}$) as $n \to \infty$.

\noindent (c) For every $r > 0$ and every 
$n \in {\mathbb N}$,
$P(\sum_{k=1}^n U_r(X_k) = 0) \geq 1 - \lambda$.}
\bigskip

   By (2.7), analogs of property (b) automatically
also hold for the dependence coefficients
$\alpha({\bf {X}},n)$ and $\rho({\bf {X}},n)$. 
Of course property (c) prevents 
the random sums in (2.13) from converging to a 
(non-degenerate) normal law under  
any choice of the $r(n)$'s and any kind of 
normalization.
Properties (a) and (c) take their most ``extreme''
forms when the parameter $\lambda$ is very small.
From (2.7) and Theorem 2.1 itself, one can see immediately that in the example described in Theorem 2.2,
$\rho(1) = \rho^*(1) = 1$ must hold.
\medskip

   Theorem 2.2 will be proved with a construction in 
Section 5.  
(Also, in Remark 5.1 at the end of that section,
a direct brief explanation of the equality 
$\rho(1) = \rho^*(1) = 1$
will be given for that particular construction.)
\medskip

{\bf Remark 2.3.}  Under its dependence 
condition (ii),
Theorem 2.1 can be extended to a weak invariance
principle (we shall not formulate that here) via the
combining of our proof (in Section 4) with an
application of the main result of Utev and 
Peligrad [2003].
Such an extension of Theorem 2.1 under its dependence condition (i) to a weak invariance principle, if that holds, would apparently require more effort, apparently involving intricate tightness arguments similar to those 
in Shao [1989] [1993].        
\bigskip

\noindent \textbf{3. Background remarks}

\medskip
For ease of reference, some elementary properties of the 
shrinking operators and corresponding expected values
are collected in the following two remarks.
(The proofs are simple and are mostly left to the reader.)
\medskip

{\bf Remark 3.1.}  For any $x\in \Rset$ and any $r,s \in [0,\infty)$ one has the following:

(a) $U_r(U_s(x))=U_{r+s}(x)$, \ \  
(the one-parameter semigroup property).

(b) $U_0(x) = x$,\ \ $U_r(-x) = -U_r(x)$, \ \ $|U_r(x)| \downarrow 0$ as $r \uparrow \infty$.

(c) $|U_r(x)|= \max\{|x|-r,0\}$, \ \ $|U_r(x)| \leq |x|$.

(d) $|U_s(x)-U_r(x)|\le|s-r|$; \ in particular, \ \ $|x-U_r(x)|\le r$.

(e) $|U_r(x)|\le |s-r|+|U_s(x)|$ (with equality
if $|x| \geq s \geq r$).

(f) If $t> 0$, then 
$|U_r(x)|\ge t$ \ if and only if \  $|x|\ge r+t$.

(g) If $\varepsilon > 0$ and $|x|\ge r+\varepsilon$, 
then $|U_r(x)|\le 2|U_{r+\varepsilon/2}(x)|$.
(That holds because $|U_r(x)| \geq \varepsilon$ and hence
$|U_r(x)| - |U_{r + \varepsilon/2}(x)| = \varepsilon/2
\leq (1/2)|U_r(x)|$.)
\medskip

{\bf Remark 3.2.}  For a given random variable $Y$
and values $r \in [0, \infty)$, one has the following: 

    (a) $E[|U_r(Y)|^2] \leq E[Y^2]$ (possibly infinity)
by Remark 3.1(c).
If\ \ \break 
$E[|U_r(Y)|^2] < \infty$ for some $r \geq 0$, 
then $E[Y^2] < \infty$ by Remark 3.1(d), and
$E[|U_r(Y)|^2] < \infty$ for all $r \geq 0$.

    (b) By a well known application of Fubini's Theorem,
followed by Remark 3.1(f),
$$ E[|U_r(Y)|^2]\
=\ 2 \int_0 ^\infty t P(|U_r(Y)| > t) dt\
=\ 2 \int_0 ^\infty t P(|Y| > t + r) dt.   \eqno (3.1)
$$

    (c) If $E[Y^2] < \infty$, then (for example by 
Remark 3.1(b)(c) and dominated convergence),
$E[|U_r(Y)|^2] \to 0$ as $r \to \infty$.
\bigskip

\noindent \textbf{4. Proof of Theorem 2.1}
\medskip

   Assume the hypotheses (and notations) in the first
paragraph of the statement of Theorem 2.1.
As was noted in Section 2, eqs.\ (2.14) and (2.15)
hold by (and are in fact equivalent to) (2.10) and (2.11),
by (3.1).
Henceforth (2.14) and (2.15) will be used freely.
In particular, from (2.14), 
$E[|U_r(X_0)|^2] < \infty$ for every
$r \in [0, \infty)$.
As in the statement of
Theorem 2.1, for each $r \in [0, \infty)$,
define the real number
$$     m_r := E[U_r(X_0)] = E[U_r(X_0) I(|X_0| > r)]
  \eqno (4.1)  $$
(where the latter equality holds by (2.8)).
Now assume that (at least) one of the weak
dependence conditions (i) or (ii) in (the third
paragraph) of Theorem 2.1 holds.
Our remaining task is to prove the last sentence
(last paragraph) of Theorem 2.1.
\medskip

   To simplify that argument, we shall convert to
an appropriate ``mean 0'' context.
For each $r \in [0, \infty)$, define 
the random sequence 
${\mathbf{Y}}^{(r)} := (Y_{k,r}, k \in {\mathbb{Z}})$
as follows:  For each $k \in \mathbb{Z}$,
$$  Y_{k,r} := U_r(X_k) - m_r.   \eqno (4.2) $$  
For each $r \in [0, \infty)$, this sequence 
${\mathbf{Y}}^{(r)}$ is strictly stationary, by the
assumption of the strict stationarity of the
sequence ${\bf X}$ in Theorem 2.1.
Further, for each $r \geq 0$ and each $n \geq 1$,
$\alpha({\bf Y}^{(r)}, n) \leq \alpha({\bf X}, n)$
and $\rho({\bf Y}^{(r)}, n) \leq \rho({\bf X}, n)$
and also 
$\rho^*({\bf Y}^{(r)}, n) \leq \rho^*({\bf X}, n)$ 
(since for each $k$ and $r$, $Y_{k,r}$ is a Borel 
function of $X_k$).
Also, for each $r \geq 0$ and each 
$k \in {\mathbb{Z}}$, 
$$ E[Y_{k,r}]=E[Y_{0,r}] = 0 \quad {\rm and} \quad
E[Y^2_{k,r}]=E[Y_{0,r}^2] = {\rm Var}\thinspace [U_r(X_0)] > 0.
\eqno (4.3) $$
Here the last inequality (>0) must hold because otherwise
one would have $U_r(X_0)=m_r$ a.s.\ by (4.1), hence 
$P(|X_0| > |m_r| + r) = 0$,
and thus $G(|m_r| + r) = 0$ (see (2.9)), which contradicts
(2.10).
\medskip

   To complete the proof of Theorem 2.1, our task
is to prove that there exists a sequence 
$(r(n), n \in {\mathbb N})$ of positive numbers satisfying
(2.12) such that
$$ 
\Bigl[\sum_{k=1}^n U_{r(n)}(X_k)\Bigl]\ -\ n\cdot m_{r(n)}\ 
=\ \sum_{k=1}^n Y_{k,r(n)}\ \Rightarrow\ N(0,1)
\quad {\rm as}\ n \to \infty.  \eqno (4.4) $$
Here the equality holds by (4.1) and (4.2).
What remains is to prove the convergence to the 
standard normal law.  
The proof will involve the random variables
$Y_{k,r}$ and will be divided into four ``steps'':
Some preliminary work will be done in Lemma 4.1 and
Step 4.2, the sequence $(r(n))$ will be constructed in
Lemma 4.3, and then the convergence to the standard 
normal law in (4.4) will be verified in
Step 4.4, thereby completing the proof of Theorem 2.1.
\medskip  
    
{\bf Lemma 4.1.}\ \ {\sl One has that
$m_r^2 = o(E[(U_r(X_0))^2])$ as $r \to \infty$,
and (hence)
$$\frac{E[Y_{0,r}^2]}{2G(r)}\ 
= \frac{E[Y_{0,r}^2]}{E[(U_r(X_0))^2]}
= \frac{Var[U_r(X_0)]}{E[(U_r(X_0))^2]}
\to\ 1\ \ {\rm as}\ r \to \infty;
\eqno (4.5) $$
and also one has that
$$ \forall \varepsilon > 0,\ \
E[Y_{0,r}^2I(|Y_{0,r}| \geq \varepsilon)]
= o(G(r))\ \ {\rm as}\ \ r \to \infty.
\eqno (4.6) $$}

   This lemma and its proof given here
are a slightly modified version, convenient for our
context, of calculations of
Housworth and Shao [2000, pp.\ 263-264].
\medskip

 {\bf Proof.}\ \ To prove (4.5) 
(and the line preceding it), first
observe that by (4.1), for any given $r \geq 0$,
\begin{eqnarray*}
|m_r| &\leq&
E[\thinspace |U_r(X_0)|
\cdot I(|X_0| > r)] \\
&\leq& \bigl(E[|U_r(X_0)|^2]\bigl)^{1/2} \cdot
\bigl(E[(I(|X_0| > r))^2]\bigl)^{1/2} \\
&=& \bigl(E[|U_r(X_0)|^2] \bigl)^{1/2} 
\cdot \bigl(P(|X_0| > r)\bigl)^{1/2}.
\end{eqnarray*}
Since $P(|X_0| > r) \to 0$ as $r \to \infty$,
one has that $m_r^2 = o(E[(U_r(X_0))^2])$
as $r \to \infty$.
And now (4.5) follows from (4.1), (4.2), (4.3), and (2.14).

\medskip
Now let us prove (4.6).
For every $r \geq 0$ and every
$\varepsilon > 0$, by Remark 3.1(f)(g) and (2.14),
one has that

\begin{eqnarray*}
E\bigl(|U_r(X_0)|^2I(|U_r(X_0)| \geq \varepsilon) \bigl)
&=&
E\bigl( |U_r(X_0)|^2I(|X_0| \geq r + \varepsilon) \bigl) \\
&\leq&
E\bigl[ |2U_{r+\varepsilon/2}(X_0)|^2I(|X_0| \geq r + \varepsilon) \bigl] \\
&\leq&  4 E[ |U_{r+\varepsilon/2}(X_0)|^2]\
=\  8G(r + \varepsilon/2).
\end{eqnarray*}
Hence by (2.11),
$$
\forall \varepsilon > 0,\ \
E\bigl( [U_r(X_0)]^2I(|U_r(X_0)| \geq \varepsilon) \bigl)
  = o(G(r))\ \ {\rm as}\ \ r \to \infty.   \eqno (4.7)
$$
Of course from Remark 3.2(c) and (4.1) (or implicitly
from above),
$m_r \to 0\ \ {\rm as}\ r \to \infty$.
If $\varepsilon > 0$, and $r\geq 0$ is such that
$|m_r| \leq \varepsilon/2$, then by (4.2),
$ \{|Y_{0,r}| \geq \varepsilon\}
\subset \{|Y_{0,r}| \leq 2|U_r(X_0)|\}$, and hence
\begin{eqnarray*}
E[Y_{0,r}^2I(|Y_{0,r}| \geq \varepsilon)]\
&\leq&
E[(2|U_r(X_0)|)^2 I(|Y_{0,r}| \geq \varepsilon)] \\
&\leq&
E[(2|U_r(X_0)|)^2 I(2|U_r(X_0)| \geq \varepsilon)] \\
&=&
4E[|U_r(X_0)|^2 I(|U_r(X_0)| \geq \varepsilon/2)]. 
\end{eqnarray*}
Hence by (4.7), eq.\ (4.6) holds.
That completes the proof of Lemma 4.1.
\medskip

   {\bf {Step 4.2}}.  Referring to the Appendix in
Section 6, applying Lemma 6.1 or Lemma 6.2
(depending on which of the dependence assumptions
(i) or (ii) in Theorem 2.1 is assumed),
let $C$ be a positive constant such that 
$$ \forall r \geq 0,\ \forall n \in {\mathbb N}, \quad
(1/C)n E[Y_{0,r}^2] \leq 
E\Bigl[ \Bigl( \sum_{k=1}^n Y_{k,r} \Bigl)^2 \Bigl]
\leq Cn E[Y_{0,r}^2].
\eqno (4.8)
$$ 

Referring to (4.3), let $N_0$ be a positive integer
such that
$$ (1/C)N_0\,E[Y_{0,0}^2] > 1.
\eqno (4.9)
$$
(Note that by (4.2) and (2.8), 
$Y_{0,0} = X_0 - m_0 = X_0 - E[X_0]$.)
\medskip

    {\bf Lemma 4.3.}\ \ {\sl For each integer $n \geq N_0$,
there exists a positive number $r(n)$ such that
$$ \Bigl\|  \sum_{k=1}^n Y_{k,r(n)} \Bigl\|_2 = 1.
\eqno (4.10)
$$
Further, if $(r(n), n \geq N_0)$ is a sequence of
positive numbers such that (4.10) holds for all
$n \geq N_0$, then
$$ r(n) \to \infty\ \ {\rm as}\ n \to \infty,
\eqno (4.11) $$
and
$$  \forall n \geq N_0, \quad
1/(Cn) \leq E[Y_{0,r(n)}^2] \leq C/n.
\eqno (4.12)
$$}

{\bf Proof.}\ \ We shall first prove the first
sentence (with eq.\ (4.10)) of Lemma 4.3.
Suppose $0 \leq r < s$. 
Then for every
$k \in {\mathbb Z}$, one has that
$|U_r(X_k) - U_s(X_k)| \leq s-r$
by Remark 3.1(d).
Hence by (4.1) and a trivial calculation,
$|m_r - m_s| \leq s-r$.
Hence by (4.2) and a simple calculation,
for any $k \in {\mathbb Z}$,
$|Y_{k,r} - Y_{k,s}| \leq 2(s-r)$.
Hence for any $n \in {\mathbb N}$,
$ |(\sum_{k=1}^n Y_{k,r}) - (\sum_{k=1}^n Y_{k,s})|
\leq 2n(s-r)$ and hence
$\|(\sum_{k=1}^n Y_{k,r}) - (\sum_{k=1}^n Y_{k,s})
\|_2 \leq 2n(s-r)$.
Hence for any $n \in {\mathbb N}$, 
$$ \Bigl| \| \sum_{k=1}^n Y_{k,r} \|_2
- \| \sum_{k=1}^n Y_{k,s} \|_2 \Bigl|
\leq 2n(s-r).
\eqno (4.13) $$
That was shown for arbitrary $0 \leq r < s$.
It follows that
for any given $n \in {\mathbb N}$, the mapping
$r \mapsto \| \sum_{k=1}^n Y_{k,r} \|_2$ for
$r \in [0, \infty)$ is (uniformly) continuous.
\medskip

    Now by Remark 3.2(c),
$E[|U_r(X_0)|^2] \to 0$ as $r \to \infty$.
Hence by the ``second half'' of (4.3),
$\|Y_{0,r}\|_2 \to 0$ as $r \to \infty$.
Hence for any given positive integer $n$,
$\|\sum_{k=1}^n Y_{k,r}\|_2 \leq
n\cdot \|Y_{0,r}\|_2 \to 0$ as $r \to \infty$.
Also, by (4.8) and (4.9),
for any $n \geq N_0$,
$E[(\sum_{k=1}^n Y_{k,0})^2] > 1$
and hence $\|\sum_{k=1}^n Y_{k,0}\|_2 > 1$.
Hence by the Intermediate Value Theorem and the
second sentence after (4.13),
for each integer $n \geq N_0$,
there exists a positive number $r(n)$ such that
(4.10) holds.
The completes the proof of the first sentence of
Lemma 4.3.
\medskip

    Now to prove the rest of Lemma 4.3, suppose
$(r(n), n \geq N_0)$ is a sequence of
positive numbers such that (4.10) holds for all
$n \geq N_0$.
By (4.8) and (4.10), one has that for any $n \geq N_0$,
$ 
(1/C)n E[Y_{0,r(n)}^2] \leq 1
\leq Cn E[Y_{0,r(n)}^2].
$
Eq.\ (4.12) follows trivially.
All that remains is to prove (4.11).
\medskip

    Suppose instead that (4.11) fails to hold;
we shall aim for a contradiction.
There exists an infinite set
$\Lambda \subset \{N_0, N_0+1, N_0 + 2, \dots \}$
and a positive number $A$ such that
$r(n) \leq A$ for all $n \in \Lambda$.
By stationarity and the second sentence after (4.13),
applied with $n=1$, the mapping
$r \mapsto \|Y_{0,r}\|_2$ for $r \in [0, \infty)$ is
continuous.
Hence by the inequality in (4.3) and the compactness
of the interval $[0,A]$ there exists a positive
number $B$ such that for all $r \in [0, A]$,
$\|Y_{0,r}\|_2 \geq B$.
In particular,
$\|Y_{0,r(n)}\|_2 \geq B$ for all $n \in \Lambda$.
Hence by (4.8),
$E[(\sum_{k=1}^n Y_{k,r(n)})^2] \geq B^2n/C$ for all
$n \in \Lambda$.
But (since $\Lambda$ is unbounded above) that
contradicts (4.10).  Hence (4.11) must
hold after all.
That completes the proof of Lemma 4.3.
\medskip

{\bf {Step 4.4.  Conclusion of proof of Theorem 2.1}}.
Applying Lemma 4.3, let $(r(n), n \geq N_0)$ be a
sequence of positive numbers that satisfies (4.10).
By Lemma 4.3, this sequence satisfies (4.11) and (4.12).
By (4.11), (4.12) and (4.5),
$\limsup_{n \to \infty} n \cdot G(r(n)) < \infty$.
Hence by (4.6) and stationarity, the Lindeberg
condition holds:
$$ \forall \varepsilon > 0,\ \
\lim_{n \to \infty}
\sum_{k=1}^n E[Y_{k, r(n)}^2I(|Y_{k, r(n)}|
\geq \varepsilon)]
= 0. \eqno (4.14) $$
Also, by (4.10) and (4.12) and stationarity
(recall the inequality in (4.3)),
$$
\liminf_{n \to \infty}
{{E\Bigl[ \Bigl(\sum_{k=1}^n Y_{k,r(n)} \Bigl)^2 \Bigl]} \over
{\sum_{k=1}^n E[(Y_{k,r(n)})^2]}}
 > 0.  \eqno (4.15) $$

    By (4.10), (4.14), (4.15), and Theorem 6.4 in the 
Appendix (Section 6), regardless of which of the
two dependence assumptions (i), (ii) in Theorem 2.1
is assumed, eq.\ (4.4) holds.  
That completes the proof of Theorem 2.1.
\medskip

\medskip
\medskip
\noindent {\bf 5. Proof of Theorem 2.2}  
\medskip
   
   As in the hypothesis of Theorem 2.2, suppose
$0 < \lambda < 1$.  
Define the number $\theta \in (0,1/4)$ by
$$ \theta := \lambda /4. \eqno (5.1) $$

    Let ${\bf V} := (V_k, k \in {\bf Z})$
be a strictly stationary Markov chain with state
space $\{1,2,3\}$, with marginal distribution given
by
$$ P(V_0 = 1) = \frac{1}{1 + 2\theta}
\quad {\rm and} \quad
P(V_0 = 2) = P(V_0 = 3) =
\frac{\theta}{1 + 2\theta},  \ \ \eqno (5.2)
$$
and with one-step transition probabilities
$p_{ij} := P(V_1=j|V_0 = i)$, for
$i,j \in \{1,2,3\}$, given by
$$ p_{11} = 1 - \theta, \quad
p_{12} = \theta, \quad {\rm and} \quad
p_{23} = p_{31} = 1.  \eqno (5.3) $$
It is easy to check that the marginal distribution in
(5.2) is the invariant distribution for the
transition probabilities given in (5.3).
It is easy to see that this Markov chain ${\bf V}$ is
irreducible and aperiodic.

To avoid frivolous technicalities, we shall henceforth
assume that for {\it every\/} $\omega \in \Omega$ and
$k \in {\mathbb Z}$, the ordered pair
$(V_k(\omega), V_{k+1}(\omega))$ is either
$(1,1)$, $(1,2)$, $(2,3)$, or $(3,1)$.
\hfil\break

    Let ${\bf Z} := (Z_k, k \in {\mathbb Z})$ be a sequence of
independent $N(0,1)$ random variables, with this
sequence ${\bf Z}$ being independent of the Markov chain 
${\bf V}$.
\hfil\break

    Define the sequence ${\bf X} := (X_k, k \in {\mathbb Z})$
as follows:  For each $k \in {\mathbb Z}$,
$$ X_k :=
\begin{cases}
           0 & \ \ \mbox{if}  \ \  V_k = 1 \\
           Z_k &  \ \  \mbox{if}  \ \ V_k = 2 \ \ \\
           -Z_{k-1} & \ \  \mbox{if} \ \  V_k = 3.  \ \   
\end{cases}
\eqno (5.4)
$$
By an elementary argument, this sequence ${\bf X}$ is 
strictly stationary.

    For events $A$ and $B$, the notation $A \doteq B$
will mean that $P(A \triangle B) = 0$, where
$\triangle$ denotes symmetric difference.
 From (5.4) and (5.1), one has that for any given
integer $k$,
$$ \{X_k = 0\} \doteq \{V_k = 1\}.  \eqno (5.5) $$

    Now $X_0$ is square-integrable and unbounded, and
hence (2.10) holds by (2.9) and (3.1).
\medskip

    The proof that this sequence ${\bf X}$ satisfies (2.11)
is simply (with one ``trivial adjustment'')
the observation in
Housworth and Shao [2000, first sentence after Theorem 1].
For the $N(0,1)$ density function $\phi(x)$,
one has that for every $\varepsilon > 0$,
$\lim_{x \to \infty} \phi(x + \varepsilon)/\phi(x) = 0$.
Now by (5.4) and (5.2), the distribution of $X_0$
is a mixture (convex combination) of the 
$N(0,1)$ distribution and the
point mass at 0; and hence it follows by symmetry and a
simple calculation that for every $\varepsilon > 0$,
$\lim_{x \to \infty}
P(|X_0| > x + \varepsilon)/P(|X_0| > x) = 0$.
It follows from a further simple integration that
(2.11) holds.
\medskip

    Now what remains is to prove properties (a), (b),
and (c) in the statement of Theorem 2.2.
We shall start with property (b).
\medskip

   {\it Proof of property (b).}\ \
Since $\rho^*({\bf Z},1) = 0$
(because of the independence of the $Z_k$'s), 
one has that
for each $n \geq 2$, by (5.4) and
Bradley [2007, v1, Theorem 6.1],
$$ \rho^*({\bf X},n) 
\leq \max\{ \rho^*({\bf V},n),\, \rho^*({\bf Z},n-1)\}
= \rho^*({\bf V},n). 
$$
Property (b) now follows from
Bradley [2007, v1, Theorem 7.15]
(applied to the Markov chain ${\bf V}$).
What remains is to prove properties (a) and (c).
\medskip

    {\it Proof of property (a).}\ \
Refer again to (2.1) and (2.4).
Let the events $A \in \sigma(X_k, k \leq 0)$ and
$B \in \sigma(X_k, k \geq 1)$ be arbitrary but fixed.
To complete the proof of property (a), it suffices to
prove that
$ |P(A \cap B) - P(A)P(B)| \leq \lambda$.
By (5.1) and the triangle inequality, it suffices to prove
$$ |P(A \cap \{X_0 = 0\} \cap B) -
P(A \cap \{X_0 = 0\})P(B)| \leq 2 \theta
\eqno (5.6) $$
and
$$ |P(A \cap \{X_0 \neq 0\} \cap B) -
P(A \cap \{X_0 \neq 0\})P(B)| \leq 2 \theta.
\eqno (5.7) $$
Between the absolute value signs in the left side
of (5.7), is the difference of two nonnegative terms
that are each trivially bounded above by
$P(X_0 \neq 0)$.
Hence the left side of (5.7) is bounded above by
$P(X_0 \neq 0)$.
By (5.5) and (5.2),
$P(X_0 \neq 0) = 2\theta/(1 + 2\theta) < 2\theta$.
Thus (5.7) holds.
To complete the proof of (a),
what remains is to prove (5.6).
\medskip

   In the proof of (5.6), for $\sigma$-fields 
${\cal A}$ and ${\cal B}$, the notation
${\cal A} \vee {\cal B}$ will mean the 
$\sigma$-field generated by ${\cal A} \cup {\cal B}$.
\medskip

    Now (see e.g.\ Bradley [2007, v1, Appendix, Remark A042])
by (5.5) and (5.2) the event $\{X_0 = 0\}$ is an
atom of the sigma field $\sigma(X_0)$.
Hence (see e.g.\ Bradley [2007, v1, Appendix, Theorem A036])
there exists an event $A^* \in \sigma(X_k, k \leq -1)$
(henceforth fixed) such that
$A \cap \{X_0 = 0\} \doteq A^* \cap \{X_0 = 0\}$.
By (5.5) and (5.4), one now has that
$A \cap \{X_0 = 0\} \doteq A^* \cap \{V_0 = 1\}$,
and also this latter event satisfies
$A^* \cap \{V_0 = 1\}
\in \sigma(V_k, k \leq 0) \vee \sigma(Z_k, k \leq -1)$.
Also by (5.4),
$B \in \sigma(V_k, k \geq 1) \vee \sigma(Z_k, k \geq 0)$.
Hence by Bradley [2007, v1, Theorem 6.2(a)],
\begin{eqnarray*}
& & [\rm{LHS\ of}\ (5.6)] \\
&=&
|P(A^* \cap \{V_0 = 0\} \cap B) -
P(A^* \cap \{V_0 = 0\})P(B)| \\  
&\leq&
\alpha\Bigl(\sigma(V_k, k \leq 0) 
\vee \sigma(Z_k, k \leq -1),
\sigma(V_k, k \geq 1) 
\vee \sigma(Z_k \geq 0) \Bigl) \\
&\leq& \alpha({\bf V},1) + \alpha({\bf Z},1)
\ =\ \alpha({\bf V},1) + 0\ =\ \alpha({\bf V},1) 
\qquad \qquad \qquad \qquad  (5.8) 
\end{eqnarray*}
(and hence as a simple consequence the second inequality
in (5.8) is really an equality).  
\medskip

    By a well known property of (strictly stationary) 
Markov chains
(see e.g.\ Bradley [2007, v1, Theorem 7.3(a)]),
$\alpha({\bf V},1) = \alpha(\sigma(V_0), \sigma(V_1))$.
The event $\{V_0 = 1\}$ is an atom of the
$\sigma$-field $\sigma(V_0)$, and by (5.2) it satisfies
$P(V_0 = 1) \geq 1 - 2\theta$.
Hence (see e.g.\
Bradley [2007, v1, Proposition 3.19(a)]),
$\alpha(\sigma(V_0), \sigma(V_1)) \leq 2\theta$.
Eq.\ (5.6) now follows from (5.8).
That completes the proof of property (a).
\medskip

    {\it Proof of property (c).}
This is the final task in the proof of Theorem 2.2.

    Suppose $n \in {\mathbb N}$ and $r \geq 0$.

Suppose $\omega \in \Omega$ is such that
$V_1(\omega) = V_n(\omega) = 1$.
Then the vector
$(V_1(\omega), V_2(\omega), \dots,
\allowbreak V_n(\omega))$
is either equal to $(1, 1, \dots, 1)$ or consists
of 1's occasionally interrupted by occurrences of
the ordered pair $(2,3)$.
If $k \in \{1, 2, \dots, n\}$ is such that
$V_k(\omega) =1$, then
$U_r(X_k(\omega)) = 0$ by (5.4).
If ($n \geq 2$ and) $k \in \{1, 2, \dots, n-1\}$
is such that
$(V_k(\omega), V_{k+1}(\omega)) = (2,3)$, then
by (5.4) and Remark 3.1(b),
\begin{eqnarray*}
U_r(X_k(\omega)) + U_r(X_{k+1}(\omega))
&=& U_r(Z_k (\omega)) + U_r(-Z_k (\omega)) \\
&=& U_r(Z_k (\omega)) - U_r(Z_k (\omega))
\ =\ 0.
\end{eqnarray*}
Hence $\sum_{k=1}^n U_r(X_k(\omega)) = 0$.

    What has been shown here is that
$\{V_1= V_n= 1\}
\subset \{\sum_{k=1}^n U_r(X_k) = 0\}$.
Now by (5.2),
$$P(V_0 = V_n = 1)
= 1 - P(\{V_0 \in \{2,3\}\} \cup \{V_n \in \{2,3\}\})
\geq 1 - 4\theta. $$
Property (c) now follows from (5.1).
That completes the proof of Theorem 2.2.
\bigskip

   {\bf Remark 5.1.}  For the (strictly stationary)
sequence ${\bf X}$ constructed above,
the equality 
$\rho({\bf X},1) = \rho^*({\bf X}, 1) = 1$,
pointed out in Section 2 in connection with Theorem 2.2, 
has the following brief explanation:  
By (5.2), (5.3), (5.4), and (5.5), the (nondegenerate)
random variables (indicator functions)
$I(\{X_{-1}=0\} \cap \{X_0 \neq 0\})$ and
$I(\{X_1 \neq 0\} \cap \{X_2 = 0\})$ are
equal a.s., 
and hence their correlation equals 1.    
\bigskip

\noindent \textbf{6. Appendix}
\medskip
    
    This section just gives convenient statements of
certain results in the literature on mixing that
are used in the proof of Theorem 2.1.
\hfil\break

    {\bf Lemma 6.1.}\ \ {\sl Suppose
$(q(1), q(2), q(3), \dots)$ is a non-increasing
sequence of nonnegative numbers such that
$q(1) < 1$ and $\sum_{n=1}^\infty q(2^n) < \infty$.
Then there exists a positive constant
$C = C(q(1), q(2), q(3), \dots)$ such that the
following holds:

    If ${\bf X} := (X_k, k \in {\mathbb Z})$ is a 
strictly stationary sequence of centered, 
square-integrable random variables such that
$\rho({\bf X},n) \leq q(n)$ for all $n \in {\mathbb N}$, then
for every $n \in {\mathbb N}$, one has that
$$
(1/C)n E[X_0^2] \leq E\Bigl[ \Bigl(
\sum_{k=1}^n X_k \Bigl)^2 \Bigl]
\leq Cn E[X_0^2].
$$}

    One reference for Lemma 6.1 is
Bradley [1988, Lemmas 2.1 and 2.3].
The upper bound (which does not require the
assumption $q(1) < 1$) was known earlier, e.g.\
from the work of Ibragimov [1975]; and in
Peligrad [1982, Lemma 3.4] it was
shown (with suitable reformulation) in a more
general form, not requiring stationarity.
\hfil\break

    {\bf Lemma 6.2.}\ \ {\sl Suppose
${\bf X} := (X_k, k \in {\mathbb Z})$ 
is a strictly stationary
sequence of centered, square-integrable random
variables such that $\rho^*(1) < 1$.
Then for every positive integer $n$,
$$
\frac{[1-\rho^*(1)]}{[1+\rho^*(1)]}\,E[X_0^2]
\leq
E\Bigl[ \Bigl(\sum_{k=1}^n X_k \Bigl)^2 \Bigl]
\leq
\frac{[1+\rho^*(1)]}{[1-\rho^*(1)]}\, E[X_0^2]. $$}

    Lemma 4.2 can be found in
Bradley [1992, Lemma 2.1] or
Bradley [2007, v1, Lemma 8.21].
In fact it follows implicitly from earlier results of
Moore [1963] involving spectral density.
\hfil\break

\medskip
    {\bf Context 6.3.}\ \ Suppose
$\xi := (\xi_{n,i}, n \in {\mathbb N},
i \in \{1, 2, \dots, k_n\})$
is a triangular array of centered,
square-integrable random variables, where
$(k_1, k_2, k_3, \dots)$ is a sequence of positive integers
such that $k_n \to \infty$ as $n \to \infty$.
For each positive integer $n$, define the nonnegative
number $\sigma_n$ by
$$\sigma_n^2 :=
E\Bigl[ \Bigl(\sum_{i=1}^{k_n} \xi_{n,i} \Bigl)^2 \Bigl]. $$
For each positive integer $m$,
define the following dependence coefficients:
\begin{eqnarray*}
\alpha(\xi,m) &:=&
\sup_{\{n: k_n \geq m+1\}}
\sup_{1 \leq u \leq k_n - m}
\alpha\Bigl(\sigma(\xi_{n,i}, 1 \leq i \leq u),\\
& & \qquad \qquad \qquad \qquad \qquad \qquad \qquad
\sigma(\xi_{n,i}, u+m \leq i \leq k_n)\Bigl); \\
\rho(\xi,m) &:=&
\sup_{\{n: k_n \geq m+1\}}
\sup_{1 \leq u \leq k_n - m}
\rho\Bigl(\sigma(\xi_{n,i}, 1 \leq i \leq u),\\
& & \qquad \qquad \qquad \qquad \qquad \qquad \qquad
\sigma(\xi_{n,i}, u+m \leq i \leq k_n)\Bigl); \\
\rho^*(\xi,m) &:=&
\sup_{\{n: k_n \geq m+1\}}\ 
\sup\ 
\rho\Bigl(\sigma(\xi_{n,i}, i \in S),\, 
\sigma(\xi_{n,i}, i \in T)\Bigl)
\end{eqnarray*}
where in the last line the inner supremum is taken
over all pairs of nonempty, disjoint sets
$S,T \subset \{1, 2, \dots, k_n\}$ such that
${\rm dist}(S,T) :=
\min_{s \in S, t \in T}|s-t| \geq m$.
\bigskip

    {\bf Theorem 6.4.}\ \ {\sl Assume that all assumptions in
the first sentence of Context 6.3 hold.
Then in the notations of Context 6.3, the following
statement holds:

Suppose that $\sum_{i=1}^{k_n} E[\xi_{n,i}^2] > 0$
for all $n$ sufficiently large,
and that
$$ \liminf_{n \to \infty}
\frac{\sigma_n^2} {\sum_{i=1}^{k_n} E[\xi_{n,i}^2]}
 > 0. \eqno (6.1) $$

    Suppose that for every $\varepsilon > 0$,
$$ \lim_{n \to \infty} \frac{1}{\sigma_n^2}
\sum_{i=1}^{k_n} E[\xi_{n,i}^2 I(|\xi_{n,i}| 
\geq \varepsilon)]
= 0 \eqno (6.2) $$
(where $I(\dots)$ denotes the indictor function).

    Suppose also that (at least) one of the following
two conditions (i), (ii) holds: \hfil\break
\noindent (i) $\sum_{m=1}^\infty \rho(\xi, 2^m) < \infty$; 
or \hfil\break
\noindent (ii) $\rho^*(\xi, m) < 1$ for some 
$m \geq 1$, 
and $\alpha(\xi, m) \to 0$ as $m \to \infty$.

    Then
$(1/\sigma_n) \sum_{i=1}^{k_n} \xi_{n,i}
\Rightarrow N(0,1)\ \ {\rm as}\ n \to \infty$.}
\bigskip

    Eq.\ (6.2) is the Lindeberg condition for
this context.
Under the mixing hypothesis (i), Theorem 6.4
is due to Utev [1990, Theorem 4.1].
Under the mixing hypothesis (ii), Theorem 6.4
is due to Peligrad [1996, Theorem 2.1]
(see also Theorem 2.2 in that paper).
\medskip

    When his paper Utev [1990] was published, Utev himself
pointed out (in a private communication to one of the
authors, R.C.B.) that eq.\ (4.3) in
his statement of his Theorem 4.1 was
intended to be (6.1) above, but had a serious
typographical error;
and he called attention to eq.\ (4.6) in his proof of
his Theorem 4.1 to show where the intended correct
version of his eq.\ (4.3) (namely (6.1) above)
was used.
\hfil\break

\medskip
\noindent {\bf References}

\medskip
\noindent R.C.\ Bradley [1988].
A central limit theorem for stationary $\rho$-mixing
sequences with infinite variance.
{\it Ann.\ Probab.\/} 16, 313-332. 
\hfil\break

\noindent R.C.\ Bradley [1992].
On the spectral density and asymptotic normality of
weakly dependent random fields.
{\it J.\ Theor.\ Probab.\/} 5, 355-373. 
\hfil\break

\noindent R.C.\ Bradley [2007].
{\it Introduction to Strong Mixing Conditions\/},
Volumes 1, 2, and 3.
Kendrick Press, Heber City, Utah. 
\hfil\break

\noindent W.\ Feller [1971].
{\it An Introduction to Probability Theory and its
Applications\/}, Vol.\ 2., 2nd edition.
Wiley, New York.
\hfil\break

\noindent H.\ F\"ollmer and A.\ Schied [2011].
{\it Stochastic Finance.\   
An Introduction in Discrete Time}, 
third revised and extended edition.
Walter de Gruyter \& Co., Berlin.
\hfil\break

\noindent H.O.\ Hirschfeld [1935].
A connection between correlation and contingency.
{\it Proc.\ Camb.\ Phil.\ Soc.\/} 31, 520-524. \hfil\break

\noindent E.\ Housworth and Q.M.\ Shao [2000].
On central limit theorems for shrunken random variables.
{\it Proc.\ Amer.\ Math.\ Soc.\/} 128, 261-267. \hfil\break

\noindent I.A.\ Ibragimov [1975].
A note on the central limit theorem for dependent
random variables.
{\it Theor.\ Probab.\ Appl.\/} 20, 135-141. \hfil\break

\noindent Z.J.\ Jurek [1977],
Limit distributions for sums of
shrunken random variables. 
In: \emph{Second Vilnius Conference
on Probability  Theory and  Mathematical Statistics.}
Abstract of Communications 3, pp. 95-96. 
\hfil\break

\noindent  Z.J.\ Jurek [1979], 
Properties of s-stable distribution functions.
{\it Bull.\ Acad.\ Polon.\ Sci.\ S\'er.\ Sci.\ Math.\/} 
27, 135-141.  
\hfil\break

\noindent Z.J.\ Jurek [1981].
Limit distributions for sums of shrunken
random variables.
{\it Dissertationes Math.\/} 185
(46 pages), PWN Warszawa. \hfil\break

\noindent Z.J.\ Jurek [1984],  
S-selfdecomposable probability measures as probability distributions of some random integrals.
{\it Limit Theorems in Probability and 
Statistics (Veszpr\'em, 1982)}, Vol.\ I, II, 
pp.\ 617-629.  
{\it Coll.\ Math.\ Societatis J\'anos Bolyai\/}, 
Vol. 36.
North Holland, Amsterdam.  
\hfil\break

\noindent Z.J.\ Jurek [1985].
Relations between the s-selfdecomposable and
selfdecomposable measures.
{\it Ann.\ Probab.\/} 13, 592-608. \hfil\break

\noindent Z.J.\ Jurek [2011], 
The random integral representation conjecture:
a quarter of a century later. 
{\it Lithuanian Math.\ J.\/} 51, 362-369.  
\hfil\break

\noindent Z.J.\ Jurek and J.D.\ Mason [1993].
{\it Operator-Limit Distributions in Probability Theory}.
Wiley, New York. \hfil\break

\noindent A.N.\ Kolmogorov and Y.A. Rozanov [1960].
On strong mixing conditions for stationary
Gaussian processes.
{\it Theor.\ Probab.\ Appl.\/} 5, 204-208. \hfil\break

\noindent M.\ Lo\`eve [1977].
{\it Probability Theory.}\ I, 4th edition.
Graduate Texts in Mathematics, Vol.\ 45.
Springer, New York.
\hfil\break

\noindent M.M.\ Meerschaert and H.P.\ Scheffler [2001].
\emph{Limit Distributions for Sums of Independent Random
Vectors: Heavy Tails in Theory and Practice\/},
Wiley Interscience, John Wiley \& Sons, New York.
\hfil\break

\noindent C.C.\ Moore [1963].  
The degree of randomness in a stationary time series.
{\it Ann.\ Math.\ Statist.\/} 34, 1253-1258. 
\hfil\break

\noindent M.\ Peligrad [1982].  
Invariance principles for mixing sequences of random
variables.
{\it Ann.\ Probab.\/} 10, 968-981.
\hfil\break

\noindent M.\ Peligrad [1996].
On the asymptotic normality of sequences of weak
dependent random variables.
{\it J.\ Theor.\ Probab.\/} 9, 703-715. 
\hfil\break

\noindent M.\ Peligrad [1998].
Maximum of partial sums and an invariance principle
for a class of weak dependent random variables.
{\it Proc.\ Amer.\ Math.\ Soc.\/} 126, 1181-1189.
\hfil\break

\noindent V.V.\ Petrov [1975].
{\it Sums of Independent Random Variables\/}.
Springer-Verlag, New York.
\hfil\break

\noindent M.\ Rosenblatt [1956].
A central limit theorem and a strong mixing condition.
{\it Proc.\ Natl.\ Acad.\ Sci.\ USA} 42, 43-47. \hfil\break

\noindent Q.M.\ Shao [1989].  
On the invariance principle for $\rho$-mixing sequences
of random variables.
{\it Chinese Ann.\ Math.\/} 10B, 427-433.
\hfil\break 

\noindent Q.M.\ Shao [1993].
On the invariance principle for stationary $\rho$-mixing
sequences with infinite variance.
{\it Chinese Ann.\ Math.\/} 14B, 27-42. \hfil\break

\noindent C.\ Stein [1972].
A bound for the error in the normal approximation to the distribution of
a sum of dependent random variables.
{\it Proceedings of the Sixth Berkeley Symposium on Probability and
Statistics\/}, vol.\ 2, pp.\ 583-602.
University of California Press, Los Angeles, 1972.
\hfil\break

\noindent S.A.\ Utev [1990].
Central limit theorem for dependent random variables.
{\it Probab. Theory and Math. Statist.\/}, 2, 519-528.
\hfil\break

\noindent S.A.\ Utev and M.\ Peligrad [2003].
Maximal inequalities and an invariance principle for
a class of weakly dependent random variables.
{\it J.\ Theor,\ Probab.\/} 16, 101-115. 
\hfil\break
\bigskip\bigskip\bigskip

\noindent Richard C.\ Bradley \hfil\break
Department of Mathematics \hfil\break
Indiana University \hfil\break
Bloomington \hfil\break
Indiana 47405 \hfil\break
USA \hfil\break
{\it E-mail\/} bradleyr@indiana.edu \hfil\break
\bigskip\bigskip

\noindent Zbigniew J.\ Jurek \hfil\break
Institute of Mathematics \hfil\break
University of Wroc\l aw \hfil\break
Pl.\ Grunwaldzki 2/4 \hfil\break
50-386 Wroc\l aw \hfil\break
Poland \hfil\break
{\it E-mail\/} zjjurek@math.uni.wroc.pl \hfil\break

\end{document}